\documentclass[reqno]{amsart}
 \usepackage{amsbsy,amssymb,amscd,amsfonts,latexsym,amstext,delarray,
 amsmath,graphicx, color, lineno}
 \usepackage[dvips]{epsfig}
\usepackage{hyperref}
\input xypic

\newcommand*\patchAmsMathEnvironmentForLineno[1]{%
  \expandafter\let\csname old#1\expandafter\endcsname\csname #1\endcsname
  \expandafter\let\csname oldend#1\expandafter\endcsname\csname end#1\endcsname
  \renewenvironment{#1}%
     {\linenomath\csname old#1\endcsname}%
     {\csname oldend#1\endcsname\endlinenomath}}%
\newcommand*\patchBothAmsMathEnvironmentsForLineno[1]{%
  \patchAmsMathEnvironmentForLineno{#1}%
  \patchAmsMathEnvironmentForLineno{#1*}}%
\AtBeginDocument{%
\patchBothAmsMathEnvironmentsForLineno{equation}%
\patchBothAmsMathEnvironmentsForLineno{align}%
\patchBothAmsMathEnvironmentsForLineno{flalign}%
\patchBothAmsMathEnvironmentsForLineno{alignat}%
\patchBothAmsMathEnvironmentsForLineno{gather}%
\patchBothAmsMathEnvironmentsForLineno{multline}%
}

\definecolor{Green}{rgb}{0,1,0}
\definecolor{Blue}{RGB}{0,0,191}
\definecolor{mathmodecolor}{RGB}{0,102,0}
\definecolor{keywordcolor}{RGB}{0,51,151}
\definecolor{sourcebackgroundcolor}{RGB}{255,247,223}
\definecolor{unixagred}{RGB}{255,0,0}
\definecolor{lightgray}{RGB}{191,191,191}
\definecolor{green}{RGB}{1,191,191}

\newtheorem{thm}{Theorem}[section]
\newtheorem{prop}[thm]{Proposition}
\newtheorem{cor}[thm]{Corollary}

\def\qqq{\,,\quad~\forall}

\def\\adsod{{\rm \adsod}}

\def\K{{\mathbb K}}

\def\N{{\mathbb N}}

\def\R{{\mathbb R}}

\def\H{{\mathbb H}}

\def\c\ads{{\mathcal \ads}}

\def\cP{{\mathcal P}}

\newcommand{\ie}{{\it i.e.\/}\ }

\newcommand{\cf}{{\it cf.\/}\ }

\def\rmax{\R_+^{\rm max}}
\def\ads{\H_\K}

\newcommand{\nil}[1]{}

\title
{Symmetrization of mono\"{\i}ds as hypergroups}

\author[Henry]{Simon Henry}

\address{Simon.~Henry: \\
\\ Paris, F-75005 France
\\ 
} \email{ henry\@@phare.normalesup.org}

\keywords{mono\"{\i}ds, hypergroups}
\subjclass[2000]{20N20,43A62}

\begin{document}

\maketitle

\begin{abstract}
We adapt the construction of the Grothendieck group associated to a commutative mono\"{\i}d to handle idempotent mono\"{\i}ds. Our construction works for a restricted class of commutative mono\"{\i}ds, it agrees with the Grothendieck group construction in many cases and yields a hypergroup which solves the universal problem for morphisms to hypergroups.
It gives the expected non-trivial hypergroup construction in the case of idempotent mono\"{\i}ds.
\end{abstract}

\tableofcontents

\section{Introduction}

The Grothendieck group construction associates an abelian group to a commutative mono\"{\i}d in a universal way.
It is used in a large variety of mathematical situations, as for example in the construction of $K$-theory groups. But as soon as the mono\"{\i}d we start from is not cancellative the Grothendieck construction looses information, and in the case of an idempotent mono\"{\i}d, \ie a mono\"{\i}d fulfilling $x+x=x$ for all $x$, this construction always yields the trivial group. It is known however, especially in the case of an idempotent mono\"{\i}d (see \ref{idem}), that there is a natural non-trivial associated hypergroup which carries more information.

We develop a notion of symmetrization which agrees with the usual one for mono\"{\i}ds such as $\N$ or $\R_+$, gives the expected non-trivial hypergroup answer in the case of idempotent additive mono\"{\i}ds such as $\rmax$, and solves the universal problem for morphisms to hypergroups.
Our main result is Theorem \ref{thm} which gives a necessary and sufficient condition on a commutative mono\"{\i}d so that the proposed symmetrization yields an hypergroup. This condition means that the set of decompositions of any element, partially ordered by refinement, is filtering, so that any two decompositions admit a common refinement. It holds for the monoid of positive elements in a commutative totally ordered group and also for the idempotent  mono\"{\i}d associated to a totally ordered set endowed with the operation $x\vee y:= {\rm max}(x,y)$. 
In both cases the construction yields the desired result.

\section{Notations}

We start by reviewing the notion of a canonical hypergroup $H$ \cite{Marty}, \cite{Krasner1}. For our applications it will be enough to consider  this particular class of  hypergroups and we shall simply call them hypergroups. We denote by $+$ the hyper-composition law in $H$. The  sum $x+y$ of two elements in $H$ is no longer a single element of $H$ but a non-empty subset of $H$. The hyper-operation on $H$ is a map
\[
+: H\times H \to \mathcal P(H)^*
\]
 taking values into the set $\cP(H)^*$ of all non-empty subsets of $H$. Thus, $\forall a,b\in H$, $a+b$ is a non-empty subset of $H$, not necessarily a singleton. One uses the notation $\forall A,B\subseteq H,~A+B:=\{\cup (a+b)|a\in A, b\in B\}$. The definition of a  (canonical) hypergroup requires that $H$ has a neutral element $0\in H$ (\ie an additive identity) and that the following axioms apply:\vspace{.05in}

 $(1)$~$x+y=y+x,\qquad\forall x,y\in H$\vspace{.05in}

$(2)$~$(x+y)+z=x+(y+z),\qquad\forall x,y,z\in H$ \vspace{.05in}

$(3)$~$0+x=x= x+0,\qquad \forall x\in H$\vspace{.05in}

$(4)$~$\forall x  \in H~  \ \exists!~y(=-x)\in H\quad {\rm s.t.}\quad 0\in x+y$\vspace{.05in}

$(5)$~$x\in y+z~\Longrightarrow~ z\in x-y.$\vspace{.05in}

 Property $(5)$ is usually called {\em reversibility}. It follows from the other $4$ axioms since 
$$
x\in y+z\implies 0\in (-x)+(y+z)\implies 0\in (-x+y)+z\implies z\in x-y
$$

\section{Universal property}

Let $B$ be a  commutative mono\"{\i}d denoted  additively  endowed with a neutral element $0$. Given a   hypergroup $H$, a map $f:B\to H$ is additive if it fulfills the following conditions:
$$ f(0)=0 $$
$$f(a+b)\in f(a)+f(b)\subset H \qqq a,b \in B$$

We say that a pair $(H,f)$ of a hypergroup $H$ and a map $f:B\to H$ is {\em universal} when  for any additive map $g:B\to K$ to a hypergroup $K$ there exists a unique morphism of hypergroups $h:H\to K$ such that $g=h\circ f$. If it exists such a pair is unique up to canonical isomorphism.

\section{Definition of $s(B)$}

Let $B$ be a mono\"{\i}d. One lets  $s(B)$ be the set $(B \times \{-1,1\}) /\sim$ , where $\sim$ is the equivalence which identifies $(a,+1) \sim (b,-1)$ for any pair $(a,b)$ such that $a+b=0$. The uniqueness of the additive inverse if it exists shows that this is an equivalence relation. One endows $s(B)$ with the hypersum $+$ given by
\begin{equation}\label{hypersum}
\begin{array}{c c c}
(a,1) + (b,1) &= & \{(a+b,1)\} \\

(a,-1)+ (b,-1) &=& \{(a+b,-1)\} \\

(a,1) + (b,-1) &=& \{(z,1) | z+b=a \} \cup \{(z,-1) | z+a=b  \} 

\end{array}
\end{equation}
One needs to show that this definition is compatible with the equivalence relation. Assume that  $b$ has an additive inverse $b'$ and let us show that $(a,1)+(b,-1) = (a,1)+(b',1)$. The required compatibility between $+$ and $\sim$ then follows. One has
\[ (a,1)+(b,-1) = \{(z,1) | z+b=a \} \cup \{(z,-1) | z+a=b  \} \]
since $b$ has an additive inverse there exists a most one $z\in B$ such that $z+b=a$ namely $z=a+b'$.
If there exists $z\in B$ such that $z+a=b$ then $a$ has an additive inverse $a'=z+b'$ and $z = a'+b$ so that $(z,-1)$ is equivalent to the other solution $(a+b',1)$, thus in the quotient one has $(a,1)+(b,-1) = (a,1)+(b',1)$. 

One has a canonical 
map $i : B \rightarrow s(B)$ which maps  $b$ to $(b,1)$, and it is an injection by construction.

To a mono\"{\i}d $B$ corresponds an ordered set $B_\leq$ associated to the preorder on $B$:
\begin{equation}\label{order}
x\leq y \iff \exists z\in B, \ x+z=y
\end{equation}
By construction, the hyperlaw \eqref{hypersum} makes sense, \ie gives a non-empty set for any pair, if and only if the order $B_\leq$ is a {\em total} order. 
\section{Theorem}

\begin{thm}\label{thm} $(i)$~Let $B$ be a mono\"{\i}d, then $s(B)$ is an hypergroup  if and only if $B$ fulfills the following property:
\begin{equation}\label{share}
\forall x,y,u,v \in B, \  \ x+y = u+v\implies \exists z \in B, \  \
\begin{cases} x+z=u, & z+v=y, 
\\
\text{or}\  \ x=u+z, & v=z+y.
\end{cases}
\end{equation}

$(ii)$~Under this hypothesis $(s(B),i)$ is the universal solution to the embedding of $B$ in a hypergroup.
\end{thm}
\proof 
$(i)$~Assume first that  $s(B)$ is an hypergroup. Let $x,y,u,v \in B$ be such that $x+y=u+v$.  We use the notation $\pm a:=(a,\pm 1) $ for the elements of $s(B)$ associated to an element $a\in B$. One then has in $s(B)$ using \eqref{hypersum}, that  $u \in (x+y)-v$, and thus using  associativity, that $u \in x+(y-v)$. 
Thus there exists  $t\in (y-v)$ such that $u \in x+t$. There are two cases:
If $t$ is an element of $B$, one has $u=x+t$ since $u \in \{x+t\}=x+t$ and $y=t+v$, using \eqref{hypersum}, since $t \in y-v$. 
Otherwise one has $t=-z$ with $z\in B$. Then $u+z=x$ since $u \in x+t=x-z$ and $z+y=v$ since $-z=t \in y-v$. Thus $B$ fulfills  condition \eqref{share}.

Let us now assume that $B$ fulfills  condition \eqref{share}. First taking $u=y$, $v=x$ we find that for any pair $x,y\in B$ there exists $z\in B$ such that $x+z=y$ or $x=y+z$. Thus $B_\leq$ is a {\em total} order and the hyperlaw \eqref{hypersum} is well defined.

 Let us show that it is associative, \ie that for $a,b,c \in s(B)$ one has $(a+b)+c = a+(b+c)$. If $a$, $b$ and $c$ have the same sign this follows from the associativity in $B$. One can thus assume that only one of the three is negative. Let us first show that $(a+b)-c = a+(b-c)$. Let $d \in (a+b)-c$ be positive, so that $d+c=a+b$. There are two cases:
 
If there exists $z\in B$ such that $b=z+c$ and $d=a+z$ then $z \in (b-c)$ and since $d=a+z$ one has $d \in a+(b-c)$.

If there exists $z\in B$ such that $c=z+b$ and $a=d+z$ then $-z \in (b-c)$ and $d \in a-z$ thus $ d \in a+(b-c)$.

Assume now that $-d \in (a+b)-c$, with $d\in B$,  \ie $d+a+b=c$. Then $-(d+a) \in (b-c)$ and $-d \in a - (d+a)$ thus $-d \in a+(b-c)$ which shows that $(a+b)-c \subset a+(b-c)$.

Let us show that $a+(b-c) \subset (a+b)-c$. The elements of $a+(b-c)$ are of three  types : the sum of $a$ and a positive element of  $(b-c)$, the positive elements of the sum of $a$ with a negative element of $(b-c)$, the negative elements of the sum of $a$ with a negative element of $(b-c)$. We treat each case separately:

Let $d=a+u$ with $u+c=b$ then $d+c=a+u+c=a+b$ \ie $d \in (a+b)-c$.

Let $d \in a-u$ with $u+b=c$ then $a=d+u$ and $a+b=d+u+b=c+d$ \ie $d \in (a+b)-c$.

Let $-d \in a-u$ with $u+b=c$ then $u=a+d$ thus $d+a+b=c$ \ie $-d \in (a+b)-c)$.

We have thus shown that  $(a+b)-c = a+(b-c)$.
It follows that $(-a+b)+c=-a+(b+c)$ since $c+(b-a) =(b+c)-a$. In the same way one gets $(a-b)+c=a+(-b+c)$ since $c+(a-b) = (c+a) -b = (a+c)-b =a+(c-b)$. We have thus shown that the hyperlaw \eqref{hypersum} is associative.

One has $(0,1)\sim (0,-1)$ and the condition $(3)$ follows. Let us check $(4)$. One has $0\in (a,1) + (b,-1) $ if and only if $a=b$. Moreover if $0\in (a,1) + (b,1) $, then $a+b=0$, so $b=a'$ is the additive inverse of $a$ and the equivalence relation identifies $(b,1)\sim (a,-1)$. Thus one gets uniqueness in all cases.

$(ii)$~Let us show that  $(s(B),i)$ is the universal solution to the embedding of $B$ in a hypergroup.
Let $H$ be an hypergroup and $f : B \to H$ be additive. Let us show that there exists a unique extension $\tilde{f}$ of $f$ to $s(B)$ such that $\tilde{f}(u+v) \subset \tilde{f}(u)+\tilde{f}(v)$.  Define $\tilde{f}$ by $\tilde{f}(a,1) = f(a)$ and $\tilde{f}(a,-1) = -f(a)$. It is well defined \ie compatible with the equivalence relation, since if 
 $u$ has an additive inverse $u'$ in $B$ one has $0 \in f(u)+f(u')$ and $f(u') = -f(u)$.  It is unique since a morphism of hypergroups fulfills $h(-x)=-h(x)$ automatically, while $\tilde{f}=f$ on $B$.

Let us check that $\tilde{f}(u+v) \subset \tilde{f}(u)+\tilde{f}(v)$.
Let $x \in u+v$, to show that $\tilde{f}(x) \in \tilde{f}(u)+\tilde{f}(v)$ one considers two cases:

If $u$ and $v$ have the same sign $u=(u',1)$, $v=(v',1)$ one has $x = (u'+v',1)$ and $\tilde{f}(x)=f(u'+v') \in f(u')+f(v')$ by additivity of $f$. (same argument if both are negative).

If $u$ and $v$ have opposite signs one can assume $u=(u',1)$, $v=(v',-1)$ and $x=(x',1)$, so that  $x'+v'=u'$.  Then $f(u') \in f(x')+f(v')$ by additivity of $f$.  Thus $\tilde{f}(u) \in \tilde{f}(x) - \tilde{f}(v) $ and one gets $\tilde{f}(x) \in \tilde{f}(u)+ \tilde{f}(v)$ using reversibility $(5)$.\endproof

In order to understand better the meaning of the condition \eqref{share} of Theorem \ref{thm}, let us reformulate it using the following notion of refinement of a decomposition into a sum $a=\sum_{i\in I}a_i$ of an element $a\in B$, where the index set $I$ is a finite totally ordered set. We say that $a=\sum_{j\in J}b_j$ is a refinement of the above decomposition if the finite totally ordered set $J$ can be decomposed as the disjoint union of intervals $J=\cup_{i\in I} J_i$ such that for each $i\in I$ one has $a_i=\sum_{j\in J_i}b_j$. This relation is symmetric and transitive and one has:
\begin{prop} \label{equiv} Let $B$ be  a commutative mono\"{\i}d. The condition \eqref{share} holds if and only if the partially ordered set of decompositions of any element of $B$ is filtering.
\end{prop}
\proof To show that the filtering condition implies condition \eqref{share}, let $x,y,u,v\in b$ such that $x+y=u+v$ and let  $a=\sum_{1\leq j\leq n}b_j$ be a common refinement of the decompositions $a=x+y=u+v$. Then there exists $k,\ell\in \{1, \ldots,n-1 \}$ such that $x=\sum_1^k b_j$, $y=\sum_{k+1}^n b_j$,  $u=\sum_1^\ell b_j$, $v=\sum_{\ell+1}^n b_j$. Using $z=\sum_{k+1}^\ell b_j$ when $k<\ell$ and $z=\sum_{\ell+1}^k b_j$ when $k>\ell$ gives the required decomposition. Conversely assume that the condition \eqref{share} holds, then let $a=\sum_{1\leq i\leq n}a_i$, $a=\sum_{1\leq j\leq m}b_j$ be two decompositions of the same element. When $n+m\leq 4$ the existence of a common refinement follows from condition \eqref{share}. We assume that we have already handled the case $n+m<N$ and consider the case $n+m=N$.  We can assume that both $n$ and $m$ are $>1$. We let 
$$
x=\sum_{1\leq i\leq n-1}a_i, \  y=a_n, \ \  u= \sum_{1\leq i\leq m-1}b_j, \   v=b_m.
$$
We apply condition \eqref{share} to $x+y=u+v$. If there exists $z\in B$ such that  $x+z=u$ and $v+z=y$, then the equality
$$
u=\sum_{1\leq i\leq n-1}a_i +z= \sum_{1\leq i\leq m-1}b_j
$$
and the induction hypothesis on $N$ give a common refinement of the form 
$$
u=\sum_{j\in J}u_j
$$
and it is enough to adjoin $v=b_m$ as the last element to obtain the required common refinement of 
$$
a=\sum_{1\leq i\leq n}a_i=\sum_{1\leq j\leq m}b_j
$$
If there exists $z\in B$ such that  $x=u+z$ and $v=y+z$, then the equality
$$
x= \sum_{1\leq i\leq m-1}b_j+z=\sum_{1\leq i\leq n-1}a_i 
$$
and the induction hypothesis on $N$ give a common refinement of the form 
$$
x=\sum_{j\in J}x_j
$$
and it is enough to adjoin $y=a_n$ as the last element to obtain the required common refinement.\endproof

Of course, in some good case, our construction produces the same group as the classical Grothendieck construction :

\begin{prop} \label{grothendieck}
Let $B$ be a mono\"{\i}d satisfying  condition \eqref{share}, then the following conditions are equivalent :

\begin{enumerate}
\item $s(B)$ is a group.
\item $s(B)$ is the group given by the Grothendieck construction (up to unique isomorphism).
\item $B$ is cancellative, i.e. $x+a=y+a \Rightarrow x=y$.
\end{enumerate}
\end{prop}

\proof
$(2) \Rightarrow (1)$ is obvious. If $s(B)$ is the group, then it is the universal solution of the problem of sending $B$ into the group, it is hence canonically isomorphic to the group of the Grothendieck construction which is also solution to this universal problem. One has $(1) \Leftrightarrow (2)$.

Assume $(1)$, if $x+a=y+a=z$ in $B$ then in $s(B)$, $x,y \in z-a$ hence if $S(B)$ is a group, one has $x=y$ which show that $(1) \Rightarrow (3)$.

Conversely, assume $(3)$, we have to show that $s(B)$ is a group, i.e. that every sum give a singleton. The sum of two positive elements or two negative elements is a singleton by definition, so we only have to consider a sum of the form $(a,1)+(b,-1) = \{(z,1) | z+b=a \} \cup \{(z,-1) | z+a=b  \}$. Consider two elements in this sum : if they are of the same sign then the cancellation condition shows they are equal, and if they are of opposite sign the cancellation condition shows that they are equivalent.

\endproof

\section{Examples}

\begin{prop}\label{usual}
Let $A$ be an abelian totally ordered group and $P$ the mono\"{\i}d of positive elements $P=\{x\in A\mid x\geq 0\}$. Then $P$ fulfills the condition of Theorem \ref{thm}, and $s(P)$ is the group $A$.
\end{prop}
\proof Let $x,y,u,v\in P$ and $x+y=u+v$. If $u\geq x$, let $z=u-x\in P$, then $y-v=z$ and one gets that there exists $z\in P$ such that  $x+z=u$ and $v+z=y$. If $u< x$ let $z=x-u\in P$, then $v-y=z$ and one gets that there exists $z\in P$ such that  $x=u+z$ and $v=y+z$. 

It remains to prove that $s(P)$ is $A$. By \ref{grothendieck}, as $P$ is cancellative, $s(P)$ is a group.  In $S(P)$ a relation $(a,1)=(b,-1)$ with $a$ and $b$ positive, means that $a+b=0$ which implies $a=b=0$ so the equivalence relation involved in the definition of $s(P)$ only identifies $(0,1)$ with $(0,-1)$. Hence there is a bijection between elements of $s(P)$ and elements of $A$ : a positive element $a \in A$ corresponds to $(a,1)$, a negative element to $(-a,-1)$, and $0$ to $(0,1)=(0,-1)$. 

The group laws on $A$ and $s(P)$ agree : When we sum same sign elements then the two additions are reduced to  addition in $P$, and in the case of elements of different signs the sum in $s(P)$ is given by $(a,1)+(b,-1)= \{(z,1) | z+b=a \} \cup \{(z,-1) | z+a=b  \} $. There are three cases :
\begin{itemize}
\item If $a-b$ is positive then the first set contains only $(a-b,1)$ and the second is empty.
\item If $a-b$ is negative then the first set is empty and the second reduces to $(b-a,-1)$.
\item If $a=b$ then $(a-b)$ is $\{0\}$ in $s(P)$.
\end{itemize}
And in all these three cases the result agrees with the sum computed in $A$.
\endproof

\begin{prop}\label{idem}
Let $S$ be a mono\"{\i}d such that $x+x=x$ for all $x\in S$. Then $S$ fulfills the condition of Theorem \ref{thm} if and only if the partial order $S_\leq$ is total.

If this holds then, then the elements of $s(S)$ are exactly :

\begin{itemize}
\item $0$,
\item The elements of the form $+s$ with $s \in S$ non zero,
\item The elements of the form $-s$ with $s \in S$ non zero.
\end{itemize} 
Moreover, the addition in the hypergroup $s(S)$ is described as follows: 
\begin{equation}\label{taurgood}
    x + y=\left\{
                 \begin{array}{ll}
                  x, & \hbox{if $|x|>|y|$ or $x=y$;} \\
                   y, & \hbox{if $|x|<|y|$ or $x=y$;} \\
                    $[$-x,x$]$ =\{(t,\pm 1) | t \leqslant |x| \}, & \hbox{if $y=-x$.}
                 \end{array}
               \right.
\end{equation}
Where $| \_ |$ denotes the map which sends $\pm a \in s(S)$ to $a \in S$.

\end{prop}
\proof We already know that the condition of Theorem \ref{thm} implies that the partial order $S_\leq$ is total. Conversely, one has, for $x,y\in S$, 
$$
x\leq y \iff x+y=y
$$
since if $x+b=y$ for some $b\in S$ one gets $x+y=x+(x+b)=(x+x) +b=y$. Thus one has $$
x\leq y \  \ \& \ \  y\leq x \implies x=y$$
and  if the partial order $S_\leq$ is total, one gets that $S$ is a totally ordered set with a minimal element $0$ and endowed with the law $x+y:= \sup(x,y)$. One needs to show that this law fulfills the condition of Theorem \ref{thm}. Thus let $x,y,u,v$ be such that  $x+y = u+v$. We need to show that there exists $z \in B$ such that  $x+z=u$ and $v+z=y$ or $x=u+z$ and $v=y+z$.
There are two cases, either the $\sup$ occurs at the same place, say $x=\sup(x,y)$, $u=\sup(u,v)$ or it occurs at different places, say $x=\sup(x,y)$, $v=\sup(u,v)$.

In the first case, assume first $y\leq v$ then take $z=v$ to get $x=u+z$ and $v=y+z$. If $v\leq y$ then take $z=y$ to get $x+z=u$ and $y=v+z$. In the second case $x=\sup(x,y)$, $v=\sup(u,v)$, take $z=x=v$. One has $x=u+z$ and $v=y+z$.

We now move to the description of $s(S)$ :
$a+b=0$ in $S$ implies that $0=a+b=a+b+b=0+b=b$ hence that $a=b=0$. So there is no other identification in the construction of $s(S)$ than $(0,1)=(0,-1)$ and the description of the elements of $s(S)$ we give in the proposition holds.

If $x$ and $y$ are both positive (or both negative), then the addition law described in the proposition reduces to the addition of $S$ and hence it coincide with the addition of $s(S)$. If we are computing a sum of the form $x-y$ with both $x$ and $y$ positives, then there are three cases :
\begin{itemize}
\item if $x>y$ then $x+z$ can never be equal to $y$ and $y+h=x$ if and only if $h=x$ hence $x-y = x$,
\item if $x<y$ then $x+z = y$ if and only if $z=y$ and $y+z$ can never be equal to $x$ hence $x-y=-y$,
\item if $x=y$ then $x+z=y$ if $z \leqslant x$ and $y=x+z$ if and only if $z \leqslant x$ hence $x-x=[-x,x]$.
\end{itemize}
Which also gives the three cases of the law described in proposition \ref{idem}.

\endproof

The next corollary illustrates the fact that in good circumstances, symmetrization commutes with dequantization:

\begin{cor}\label{exidem} Let $S$ be the mono\"{\i}d $\rmax$, then the hypergroup $s(S)=s(\rmax)$ is the dequantization $\R^\flat$ (\cf \cite{CC5}, \cite{viro}) of the 
symmetrization $s(\R_+)=\R$ of the semifield of positive real numbers (with usual operations).
\end{cor}


\begin{thebibliography}{99}







\bibitem{CC5}  A.~Connes, C.~Consani,  {\em The universal thickening  of the field of real numbers}, arXiv:0420079v1 [mathNT] (2012).


\bibitem{Krasner1} M.~Krasner, {\em A class of hyperrings and hyperfields}. Internat. J. Math. Math. Sci. 6 (1983), no. 2, 307--311.


\bibitem{Marty} F.~Marty, {\em Sur une g\'en\'eralisation de la notion de groupe} in Huiti\`eme Congres des Math\'ematiciens, Stockholm 1934, 45--59.


\bibitem{viro} O.~Viro {\em Hyperfields for tropical geometry I, hyperfields and dequantization}
arXiv:1006.3034v2.










\end{thebibliography}
\end{document}